\documentclass{article}
\usepackage{amsmath}
\usepackage{amssymb}
 
\newtheorem{thm}{Theorem}[section]
\newtheorem{pro}[thm]{Proposition}
\newtheorem{cor}[thm]{Corollary}
\newtheorem{lem}[thm]{Lemma}

\newtheorem{rem}[thm]{Remark}
\newtheorem{defn}[thm]{Definition}
\newtheorem{exam}[thm]{Example}

\newcommand{\spec}{{\rm Spec}}
\newcommand{\gr}{{\rm grade}}
\newcommand{\de}{{\rm depth}}

\newcommand{\supp}{{\rm Supp}}
\newcommand{\hgt}{{\rm ht}}
\newcommand{\ann}{{\rm Ann}}
\newcommand{\pd}{{\rm pd}}

\newcommand{\smin}{{\rm Min}}

\newcommand{\Ext}{{\rm Ext}}
\newcommand{\Tor}{{\rm Tor}}
\begin{document}
  
\title{Tensor products and direct limits of almost Cohen-Macaulay modules}

\author{Cristodor Ionescu and Samaneh Tabejamaat}

\maketitle
\date{}
\begin{abstract}
  We investigate  the almost Cohen-Macaulay property and the Serre-type  condition $(C_n),\ n\in\mathbb{N},$ for noetherian algebras and  modules.
More precisely, we find permanence properties of these conditions with respect to tensor products and direct limits.
\end{abstract}

\section{Introduction}

All rings considered will be commutative, with unit and noetherian. All modules are supposed to be finitely generated.
\par Almost Cohen-Macaulay rings appeared from a flaw in Matsumura's book \cite{Mat} and were first studied by Han \cite{Han} 
and afterwards by Kang \cite{K1}, who introduced the notion of almost Cohen-Macaulay modules. 
The first author considered the condition $(C_n),$ where $n$ is a natural number \cite{I},  inspired by the well-known condition 
$(S_n)$ 
of Serre and characterized almost Cohen-Macaulay rings using this notion. The notion of a module satisfying the  condition $(C_n)$ 
was defined 
and studied by the second author and A. Mafi \cite{MT2}. 
\par We study the behaviour of the condition $(C_n), \ n\in\mathbb{N}$  and 
of almost Cohen-Macaulayness with respect to tensor product of $A$-modules and of $A$-algebras.
  We show that a ring having a module of finite projective dimension satisfying the condition $(C_n),$ 
must satisfy itself the condition $(C_n)$, a property that was already proved for the condition $(S_n).$ We investigate also the
 behaviour of the condition
$(C_n)$ 
and of almost Cohen-Macaulayness to direct limits.

\section{Tensor products of almost Cohen-Macaulay algebras}
We  begin by recalling the notions and basic facts that will be needed in the paper.

\begin{defn}\label{acmdef} {\rm \cite[Definition 1.2]{K1}}
Let $A$ be a ring and $M$ be an $A$-module. We say that $M$ is an almost Cohen-Macaulay $A$-module if  $P\in\supp(M)$ we have $\de(P,M)=\de_{A_P}(M_P).$ 
$A$ is called an almost Cohen-Macaulay ring if it is an almost Cohen-Macaulay $A$-module.
\end{defn}

\begin{rem}\label{acmmod}
Let $(A,m)$ be a noetherian local ring and $M$ a finitely generated $A$-module. Then it follows at once by {\rm \cite[Lemma 1.5 and Lemma 2.4]{K1}} 
that $M$ is almost Cohen-Macaulay iff $\dim(M)\leq\de(M)+1.$ 
\end{rem}

\begin{defn}\label{cn} {\rm \cite[Definition 2.1]{MT2}}
Let $A$ be a ring, $n$ a natural number and $M$ an $A-$module. We say that $M$ satisfies the condition $(C_n)$ if for every $P\in\supp(M)$ we have 
$\de(M_P)\geq \min(n,\hgt_M(P))-1.$ If $A$ satifies the condition $(C_n)$ as an $A$-module, we say that the ring $A$
 satisfies $(C_n).$
\end{defn}

\begin{rem}\label{sn}
 Recall that an $A-$module $M$ is said to satisfy  Serre's condition $(S_n)$ if for every $P\in\supp(M)$ we have $\de(M_P)\geq \min(n,\hgt_M(P)).$ Hence
 obviously if $M$ satisfies the condition $(S_n),$ it satisfies also the condition $(C_n).$
\end{rem}

\begin{lem}\label{CnSn}
 Let $n\in\mathbb{N}$ and $M$ be an $A$-module. If $M$ satisfies the condition $(S_n),$ then it satisfies the condition $(C_{n+1}).$ 
\end{lem}

\par\noindent \textit{Proof:} Let $P\in\supp(M).$ Then by the condition $(S_n)$ we have $\de(M_P)\geq \min(\hgt_M(P),n).$
 Suppose first that $n<\hgt_M(P).$ Then $\min(\hgt_M(P),n)=n,$   hence $\de(M_P)\geq n.$ But $n+1\leq\hgt_M(P),$ 
hence $\min(n+1,\hgt_M(P))=n+1.$ 
Then $\min(n+1,\hgt_M(P))-1=n+1-1=n.$ It follows that $\de(M_P)\geq n=\min(n+1,\hgt_M(P))-1.$
\par Suppose now that $n\geq \hgt_M(P).$ Then $\min(n,\hgt_M(P))=\hgt_M(P)$, hence by $(S_n)$ we have $\de(M_P)=\hgt_M(P).$ 
But $n+1>\hgt_M(P),$ hence $\min(n+1,\hgt_M(P))-1= \hgt_M(P)-1<\hgt_M(P)=\de(M_P).$

\begin{lem}\label{fibresacm}
Let $u:A\to B$ be a flat morphism of noetherian rings. If $B$ is almost Cohen-Macaulay, then all the fibers of $u$ are almost Cohen-Macaulay.
\end{lem}

\par\noindent \textit{Proof:}  Let $P\in\spec(A)$ and let $QB_P/PB_P\in\spec(B_P/PB_P).$ Then the morphism $A_P\to B_Q$ is flat and local
and by \cite[Proposition 2.2,a)]{I} it follows that $B_Q/PB_Q=(B_P/PB_P)_{QB_P/PB_P}$ is almost Cohen-Macaulay. From \cite[Lemma 2.6]{K1}
we obtain that $B_P/PB_P$ is almost Cohen-Macaulay.

\begin{pro}\label{tens1}
Let $k$ be a field, $A$ and $B$ be $k-$algebras such that $A\otimes_kB$ is noetherian and $n\in\mathbb{N}.$ Then:
\par i) If $A\otimes_kB$ satisfies the condition $(C_n),$ then $A$ and $B$ satisfy the condition $(C_n);$
\par ii) If $A$ and $B$ satisfy the condition $(C_n)$ and one of them satisfies the condition $(S_n),$ then $A\otimes_kB$
 satisfies the condition $(C_n);$  
\par iii) If $A$ and $B$ satisfy the condition $(S_n)$, then $A\otimes_kB$ satisfies the condition $(C_{n+1}).$
\end{pro}

\par\noindent \textit{Proof:}  i) Follows from \cite[Proposition 3.12]{I}.
\par\noindent ii) Follows from \cite[Proposition 3.13]{I}.
\par\noindent iii) Follows from \cite[Theorem 6,b)]{TY} and Lemma \ref{CnSn}.

\begin{cor}\label{tens2}
 Let $k$ be a field, $A$ and $B$ be $k-$algebras such that $A\otimes_kB$ is noetherian. Then:
 \par i) If $A\otimes_kB$ is almost Cohen-Macaulay, then $A$ and $B$ are almost Cohen-Macaulay;
 \par ii) If $A$ and $B$ are almost Cohen-Macaulay and moreover one of them is Cohen-Macaulay, then $A\otimes_kB$ is almost 
Cohen-Macaulay.  
\end{cor}

\par\noindent \textit{Proof:}  Follows from \ref{tens1} and \cite[Theorem 3.3]{I}.

\begin{exam}\label{acm}
 Let $A=k[x^4,x^5,xy,y]_{(x^4,x^5,xy,y)}.$ Then $A$ is a noetherian local domain of dimension 2, 
 hence by 
remark \ref{acmmod} it is
  almost Cohen-Macaulay. Since $x^5y=x^4xy\in x^4A$ and $x^5(xy)^3=x^8y^3\in x^4A,$ it follows that $\de(A)=1.$ But $A\otimes_kA$ 
  is not almost Cohen-Macaulay, because by {\rm \cite[Lemma 2]{F}} we get $\dim(A\otimes_kA)=4$ and $\de(A\otimes_kA)=2.$ 
  The first example with this property was given by Taba\^a  {\rm \cite[Exemple]{T}}.
\end{exam}

\begin{pro}\label{tensprodmod}
Let $u:A\to B$ be a morphism of noetherian rings, $M$ a finitely generated $A$-module, $N$ a finitely generated $B$-module and 
$n$ a natural number.
Suppose that $N$ is a flat $A$-module. We consider the structure of $M\otimes_AN$ as a $B$-module. Then:
\par i) If $M\otimes_AN$ satisfies the condition $(C_n)$, then $M$ satisfies the condition $(C_n);$
\par ii) If $M$ and $N_P/PN_P$ satisfy the condition $(S_n),$ for any $P\in\spec(A),$ then $M\otimes_AN$ satisfies the condition $(S_n);$
\par iii) If $M$ satisfies the condition $(S_n)$ and $N_P/PN_P$ satisfies the condition $(C_n)$ for any $P\in\supp(M),$ 
then $M\otimes_AN$ satisfies the condition $(C_n);$
\par iv)  If $M$ satisfies the condition $(C_n)$ and $N_P/PN_P$ satisfies the condition $(S_n)$ for any $P\in\supp(M),$ then $M\otimes_AN$
satisfies the condition $(C_n).;$
\par v) If $M$ and $N_P/PN_P$ satisfy the condition $(S_n),$ for any $P\in\supp(M),$ then $M\otimes_AN$ satisfies the condition $(C_{n+1}).$
\end{pro}

\par\noindent \textit{Proof:}  i) Let $P\in\supp(M), Q\in\smin(PB).$ By \cite[Proposition 1.2.16 and Theorem A.11]{BH}, the flatness of $u$ implies that 
$$\de_{B_Q}(M_P\otimes_{A_P}N_Q)=\de_{A_P}(M_P)+\de_{B_Q}(N_Q/PN_Q)$$
and 
$$\dim_{B_Q}(M_P\otimes_{A_P}N_Q)=\dim_{A_P}(M_P)+\dim_{B_Q}(N_Q/PN_Q).$$
But $\dim_{B_Q}(N_Q/PN_Q)=0$ and $M_P\otimes_{A_P}N_Q$ satisfies the condition $(C_n),$ hence 
$$\de_{A_P}(M_P)\geq\min(n,\dim_{A_P}(M_P))-1,$$ that is $M$ satisfies the condition $(C_n).$ 
\par\noindent ii) Let $Q\in\supp_B(M\otimes_AN)$ and $P=Q\cap A.$ As above we have 
$$\de_{B_Q}(M_P\otimes_{A_P}N_Q)\geq\min(n,\dim_{A_P}(M_P))+\min(n,\dim_{B_Q}(N_Q/PN_Q))\geq$$
$$\geq\min(n,\dim_{B_Q}(M_P\otimes_{A_P}N_Q)).$$
\par\noindent iii) and iv) The proof is similar to the proof of assertion i).
\par\noindent v) Follow from ii) and lemma \ref{CnSn}.

\begin{pro}\label{tensprodmod1}
Let $u:A\to B$ be a morphism of noetherian rings, $M$ a finitely generated $A$-module, $N$ a finitely generated $B$-module and $n$ 
a natural number. 
Suppose that $N$ is a flat $A$-module. We consider the structure of $M\otimes_AN$ as a $B$-module. Then:
\par i) If $M\otimes_AN$ is almost Cohen-Macaulay, then $M$ and $N_P/PN_P$ are almost Cohen-Macaulay for any $P\in\supp(M);$
\par ii) If $M$ is almost Cohen-Macaulay and $N_P/PN_P$ is Cohen-Macaulay for any $P\in\supp(M),$ then 
$M\otimes_AN$ is almost Cohen-Macaulay;
\par iii) If $M$ is Cohen-Macaulay and $N_P/PN_P$ is almost Cohen-Macaulay for any $P\in\supp(M),$ then $M\otimes_AN$
 is almost Cohen-Macaulay.
\end{pro}

\par\noindent \textit{Proof:}  i) From \cite[Theorem 3.3]{I} and Proposition \ref{tensprodmod} i), it follows at once that $M$ is almost Cohen-Macaulay. 
Let $P\in\supp_A(M).$ We have  
$$\supp_{B_P/PB_P}(N_P/PN_P)=$$
$$=\{QB_P/PB_P\ \vert\ Q\in\supp_B(N)\cap V(PB), Q\cap(A\setminus P)=\emptyset\ \}.$$ 
Hence, let $QB_P/PB_P\in\supp_{B_P/PB_P}(N_P/PN_P).$ Since $M\otimes_AN$ is an almost Cohen-Macaulay $B-$module, by \cite[Lemma 2.6]{K1}
we get that $M_P\otimes_{A_P}N_Q$ is an almost Cohen-Macaulay $B_Q-$module, as $Q\in\supp(M\otimes_AN).$ Since $u$ is flat, $M_P$ is 
finitely generated and $N_Q$ is finitely generated $B_Q$-module and flat $A_P$-module, by \cite[Proposition 2.7]{MT} it follows that $N_Q/PN_Q$ 
is an almost Cohen-Macaulay $B_Q/PB_Q$-module. Now by \cite[Lemma 2.6]{K1} we obtain that $N_P/PN_P$ is an almost Cohen-Macaulay $B_P/PB_P$-module.
\par\noindent ii) and iii) Follows from \ref{tensprodmod} and \cite[Theorem 3.3]{I}.

\begin{lem}\label{l1}
Let $R$ be a commutative ring, $A$ and $B$ be $R$-algebras such that $A\otimes_RB$ is noetherian.  Let $P\in\spec(A\otimes_RB)$ and set $p:=P\cap A, q=P\cap B, r:=P\cap R.$
Assume that $A_p$ is flat over $R_r.$ Then:
\par i) If $(A\otimes_RB)_P$ is almost Cohen-Macaulay, then $B_q$ and $A_p/rA_p$ are almost Cohen-Macaulay;
\par ii) If $A_p/rA_p$ and $B_q$ are almost Cohen-Macaulay and one of them is Cohen-Macaulay, then $(A\otimes_RB)_P$ is almost Cohen-Macaulay.
\end{lem}

\par\noindent \textit{Proof:}  Follows at once from \cite[Corollary 2.8]{BCLM}.


\begin{cor}\label{l2}
Let $R$ be a commutative ring, $A$ and $B$ be $R$-algebras such that $A\otimes_RB$ is noetherian.  Let $P\in\spec(A\otimes_RB)$ and set $p:=P\cap A, q=P\cap B, r:=P\cap R.$ Assume 
that $A_p$ and $B_q$ are flat over $R_r.$ Then:
\par i) If $(A\otimes_RB)_P$ is almost Cohen-Macaulay, then $A_p$ and $B_q$ are almost Cohen-Macaulay;
\par\ ii) If $A_p$ and $B_q$ are almost Cohen-Macaulay and one of them is Cohen-Macaulay, then $(A\otimes_RB)_P$ is almost Cohen-Macaulay.
\end{cor}

\par\noindent \textit{Proof:}  i) Follows from \ref{l1} and the flatness of $A_p$ and $B_q$ over $R_r.$
\par\noindent ii) Since $R_r\to A_p$ is flat, by \cite[Proposition 2.2]{I} it follows that $A_p/rA_p$ is almost Cohen-Macaulay. 
Note that in case $A_p$ is Cohen-Macaulay, $A_p/rA_p$ is Cohen-Macaulay too. Now the assertion follows from \ref{l1}, ii). 

\begin{exam}\label{exemplu2}
In the previous corollary, the assumption in ii) that one of the rings is Cohen-Macaulay is necessary. Indeed, 
consider again the ring in example \ref{acm}, that is $A=B=k[x^4,x^5,xy,y]_{(x^4,x^5,xy,y)}$ and let $R=k[x^4].$ 
By \cite[Proposition 2.1, 2)]{BCLM}, there is $P\in\spec(A\otimes_RB)$  such that $P\cap A$ is the maximal ideal of $A.$  Then by 
\cite[Corollary 2.6]{BCLM} it follows that $(A\otimes_RB)_P$ is
not almost Cohen-Macaulay. As we saw in \ref{acm} the ring $A_p=B_q$ is  almost Cohen-Macaulay.
\end{exam}

\begin{cor}\label{l3}
Let $R$ be a commutative ring, $A$ and $B$ be flat $R$-algebras such that $A\otimes_RB$ is noetherian.  Then:
\par i) If $A\otimes_RB$ is almost Cohen-Macaulay, then  $A_p$ and $B_q$ are almost Cohen-Macaulay for any $p\in\spec(A)$ and $q\in\spec(B)$ such 
that  $p\cap R=q\cap R$;
\par ii) If for any $p\in\spec(A)$ and $q\in\spec(B)$ such that 
$p\cap R=q\cap R, A_p$ and $B_q$ are almost Cohen-Macaulay and one of them is Cohen-Macaulay, then $A\otimes_RB$ is almost Cohen-Macaulay. 
\end{cor}

\par\noindent \textit{Proof:}  i) By \cite[Proposition 2.1]{BCLM} there exists $P\in\spec(A\otimes_RB)$ such that $P\cap A=p, P\cap B=q.$ 
Now apply \ref{l2}, i).
\par\noindent ii) Follows from \ref{l2}, ii).

\section{Tensor products of almost Cohen-Macaulay modules}

Recall  the following definition(cf. \cite[Definition 1.1]{YKH}):
\begin{defn}\label{defgr} If $M$ and $N$ are finitely generated non-zero $A$-modules, the grade of $M$ with respect to $N,$ is 
defined by  $$\gr(M,N)=\inf\{i\,\vert\, \Ext^i_A(M,N)\neq 0\}.$$
\end{defn}

\begin{rem}\label{fox} By {\rm \cite[Proposition 1.2, h)]{Fo}} we see that 
$\gr(M,N)$ is the length of a maximal $N$-regular sequence in $\ann(M).$
\end{rem}

\begin{rem}\label{perfrem} By {\rm \cite[Proposition 1.2.10, a)]{BH}} we get 
$$\gr(M,N)=\inf\{(\de(N_P)\, \vert\, P\in\supp(M)\}=$$
$$=\inf\{(\de(N_P)\, \vert\,  P\in\supp(M)\cap\supp(N)\}.$$
\end{rem}

\begin{defn}\label{perf}
If $M$ and $N$ are finitely generated $A$-modules and $\pd(M)<\infty,$ we say that $M$ is $N$-perfect if $\gr(M,N)=\pd(M).$
\end{defn}

We will use several times the following fact:
\begin{rem}\label{inters}{\rm (Intersection Theorem cf. \cite[Theorem 8.4.4]{ST})}
If $\pd_A(M)<\infty,$ then $\dim(N)\leq\pd(M)+\dim(M\otimes_AN).$
\end{rem}

\begin{pro}\label{sam1}
Let $A$ be a noetherian local ring, $M$ and $N$ be finitely generated non-zero $A$-modules, $n\in \mathbb{N}.$ Assume that:
\par a) $M$ is $N-$perfect and  $\pd(M):=p\leq n;$
\par b) $N$ satisfies the condition $(C_n);$
\par c) $\Tor^A_i(M,N)=0,\ \forall i>0.$
\par\noindent Then $M\otimes_AN$ satisfies the condition $(C_{n-p}).$
\end{pro}

\par\noindent \textit{Proof:}  Let $P\in\supp(M\otimes_AN).$ We have 
$$p=\gr(M,N)\leq\gr(M_P,N_P)\leq\pd(M_P)\leq\pd(M)=p,$$
 hence $\pd(M_P)=p.$ By \cite[Theorem 1.2]{AU} we obtain 
$$\de(M\otimes_AN)_P=\de(M_P\otimes_{A_P}N_P)=\de(N_P)-\pd(M_P)=$$ 
$$=\de(N_P)-p\geq\min(\hgt_NP,n)-1-p=\min(\hgt_NP-p,n-p)-1$$
and it is enough to show that $\min(\hgt_NP-p,n-p)\geq\min(\hgt_{M\otimes_AN}P,n-p).$ But since $\gr(M_P,N_P)=\gr((M\otimes_AN)_P,N_P)$,
by  \cite[Theorem 2.1]{YKH} we have $p\leq \hgt_N(P)-\hgt_{M\otimes_AN}(P)$ and this concludes the proof.

\begin{cor}\label{perf3}
Let $n\in\mathbb{N}$, $A$ be a noetherian local ring satisfying the condition $(C_n)$ and let $M$ be  
a perfect $A$-module such that $\pd(M)=p\leq n.$
Then $M$ satisfies the condition $(C_{n-p}).$
\end{cor}

\begin{pro}\label{perf4}
Let $A$ be a noetherian local ring and let $M$ and $N$ be finitely generated non-zero $A$-modules such that
$\Tor^A_i(M,N)=0,\ \forall i>0.$ If $\pd(M)<\infty$ and  $M\otimes_AN$ is almost Cohen-Macaulay, 
then $N$ is almost Cohen-Macaulay.
\end{pro}

\par\noindent \textit{Proof:}  Applying remarks \ref{acmmod} and \ref{inters} and \cite[Theorem 1.2]{AU} we get
 $$\dim(N)\leq\pd(M)+\dim(M\otimes_AN)\leq $$
$$\leq\pd(M)+\de(M\otimes_AN)+1=\de(N)+1,$$
hence applying again remark \ref{acmmod} it follows that $N$ is almost Cohen-Macaulay.

\begin{cor}\label{acmnice} {\rm(see \cite[Corollary 8.4.5]{ST})}
Let $A$ be a  local ring having an almost Cohen-Macaulay module of finite projective dimension. Then $A$ is an almost Cohen-Macaulay ring.
\end{cor}

We can prove a more general form of \ref{perf4}:

\begin{pro}\label{perf5}
 Let $A$ be a noetherian local ring and $n\in\mathbb{N}.$ Let also $M$ and $N$ be finitely generated non-zero $A$-modules such that
$\Tor^A_i(M,N)=0,\ \forall i>0.$ If $\pd(M)<\infty$ and  $M\otimes_AN$ satisfies the condition $(C_n),$ then $N$ satisfies the condition $(C_n).$ 
\end{pro}

\par\noindent \textit{Proof:}  Let $P\in\supp_A(N).$ Assume first that $P\in\supp_A(M ).$  Then $P\in\supp(M\otimes_AN)$ and we have two cases:
\par\noindent a) $\de(M\otimes_AN)_P<n-1.$ Then $(M\otimes_AN)_P$ is almost Cohen-Macaulay and by Proposition \ref{perf4} it follows
that $N_P$ is almost Cohen-Macaulay too. 
Hence $\de(N_P)\geq\min(n, \dim(N_P))-1$.
\par\noindent b) $\de(M\otimes_AN)_P\geq n-1.$ Then  from\cite[Theorem 1.2]{AU} we obtain:
$$\de(N_P)=\de(M\otimes_AN)_P+\pd(M_P)\geq$$
$$\geq n-1\geq\min(n-1,\hgt_N(P)-1)=\min(n,\hgt_N(P))-1.$$
Assume now that $P\notin\supp_A(M).$ Let $Q\in\smin(\ann_A(M)+P).$ By \cite[Theorem 2.1]{YKH} we have 
$$\de(N_P)\geq\gr(A_Q/PA_Q,N_Q)\geq\de(N_Q)-\dim(A_Q/PA_Q).$$ 
By remark \ref{inters} we  have 
$$\dim(A_Q/PA_Q)\leq\pd(M_Q)+\dim(M_Q/PM_Q)=\pd(M_Q),$$
hence by \cite[Theorem 1.2]{AU} we obtain
$$\de(N_P)\geq\de(N_Q)-\pd(M_Q)=\de(M_Q\otimes_{A_Q}N_Q).$$
If $\de(M_Q\otimes_{A_Q}N_Q)\geq n-1$
we have $$\de(N_P)\geq n-1\geq\min(n-1,\hgt_N(P)-1)=\min(n,\hgt_N(P))-1,$$
as required.
\par\noindent In the case $\de(M_Q\otimes_{A_Q}N_Q)<n-1$, by definition $M_Q\otimes_{A_Q}N_Q$
is almost Cohen-Macaulay,  hence from Proposition \ref{perf4} it follows
that $N_Q$ is almost Cohen-Macaulay. But then $N_P$ is almost Cohen-Macaulay too and
consequently $$\de(N_P)\geq\min(n,\hgt_N(P))-1.$$


\begin{cor}\label{acmnice2}
Let $A$ be a  local ring having a module of finite projective dimension satisfying the condition $(C_n)$. Then $A$ satisfies the condition $(C_n)$.
\end{cor}

\begin{rem}\label{yossn}
 A similar result for the condition $(S_n)$ was proved in {\rm \cite[Proposition 4.1]{Y}}.
\end{rem}

\section{Direct limits of almost Cohen-Macaulay rings} 

\begin{pro}\label{limclas}
 Let $(A_{i},f_{ij})_{i,j\in\Lambda}$ be a direct system of noetherian rings, $n\in\mathbb{N}$ and let $A:=\mathop{\varinjlim}\limits_{i\in \Lambda} A_i$. Assume that:
\par a) the ring $A$ is noetherian;
\par b) the morphism $f_{ij}$ is flat for any $i\leq j;$
\par c) the ring $A_i$ satisfies the condition $(C_n)$ for any $i\in \Lambda.$
\par\noindent Then $A$ satisfies the condition $(C_n).$
\end{pro}

\par\noindent \textit{Proof:}  Let $P\in \spec(A), B:=A_P, k:=B/PB$ and for any $i\in\Lambda,$ put  $P_i:=P\cap A_i, B_i:=(A_i)_{P_i}$ and
$ k_i:=B_i/P_iB_i.$ 
There exists $i_0\in\Lambda$
such that $P=P_iB$ for any $i\geq i_0,$ hence $k=B\otimes_{B_i}k_i.$ But the morphism $B_i\to B$ is flat, hence 
$\dim(B_i)=\dim(B)$ and $\de(B_i)=\de(B)$ for any $i\geq i_0.$ The assertion follows from \ref{acmmod}.

\begin{cor}\label{limflat}
 Let $(A_{i},f_{ij})_{i,j\in\Lambda}$ be a direct system of noetherian rings and let $A:=\mathop{\varinjlim}\limits_{i\in \Lambda} A_i$.
 Assume that:
\par a) the ring $A$ is noetherian;
\par b) the morphism $f_{ij}$ is flat for any $i\leq j;$
\par c) the ring $A_i$ is almost Cohen-Macaulay for any $i\in \Lambda.$
\par\noindent Then $A$ is almost Cohen-Macaulay.
\end{cor}

\par\noindent \textit{Proof:}  Let us give a direct proof, slightly different to the one coming out at once from \ref{limclas}. 
We will apply \cite[Corollary 2.3]{CTT}.
Let $P\subseteq Q$ be two prime ideals in $A$ and for any $i\in\Lambda$ let $P_i:=P\cap A_i, Q_i:=Q\cap A_i, B_i:=(A_i)_{P_i}, 
C_i:=(A_i)_{Q_i}.$ Let also $B:=A_P$ and $C:=A_Q.$ Then, as in the proof of 
\ref{limclas}, there exists $i_0\in\Lambda$ such that  $P=P_iB$ and $Q=Q_iB$ for any $i\geq i_0.$ Since $A_i$ is almost Cohen-Macaulay 
it 
follows by \cite[Corollary 2.3]{CTT} that $\de(B_i)\leq \de(C_i)$ and consequently, by flatness, $\de(B)=\de(B_i)\leq\de(C_i)=\de(C).$ 
Now we apply again \cite[Corollary 2.3]{CTT}.

\begin{pro}\label{limite}
Let $((A_{i},m_i,k_i),f_{ij})_{i,j\in\Lambda}$ be a direct system of noetherian local rings and let $A:=\mathop{\varinjlim}\limits_{i\in \Lambda} A_i$.
Assume that:
\par a) $m_iA_j=m_j,\ \forall i\leq j;$
\par b) for any $i\in \Lambda$ and for any ideal $I$ of $A_i$ we have $IA\cap A_i=I;$
\par c) $A_i$ is an almost Cohen-Macaulay ring for any $i\in \Lambda.$
\par\noindent Then $A$ is almost Cohen-Macaulay.
 \end{pro}

\par\noindent \textit{Proof:}  By \cite{Og} and \cite[$(0_{III}, 10.3.1.3)$]{EGA},  the ring $A$ is a local noetherian ring, with maximal ideal
$m=m_iA, \forall i\in\Lambda.$  Let $i_0\in \Lambda$ be such that 
$\dim(A_j)=\dim(A),\ \forall j\geq i_0$(see \cite[Lemma 3.10]{ADT}) and let $x_1,\ldots,x_s\in A_{i_0}$ be a system of parameters. 
Set $Q_j=(x_1,\ldots,x_s)A_j, \forall j\geq i_0$ and $Q=(x_1,\ldots,x_s)A.$ 
Then $$m=m_{i_0}A=\sqrt{Q_{i_0}}A\subseteq\sqrt{Q_{i_0}A}\subseteq m,$$ whence $x_1,\ldots,x_s$ is a system of parameters in $A.$ 
The same argument shows that 
$x_1,\ldots,x_s$ is a system of parameters in $A_j,\ \forall j\geq i_0.$
If $A_i$ is CM for all $i\in \Lambda,$ it follows by  \cite[Corollary 3.7]{ADT} that $A$ is Cohen-Macaulay and we are done. 
Assume not. Since $A_j$ is almost CM, by \cite[Theorem 1.7]{K1}, there exists $l\in\{1,\ldots,s\}$ such that
$\{x_1,\ldots,\widehat{x_l},\ldots,x_s\}$ is a regular sequence. We assume $l=1,$ that is $\{x_2,\ldots,x_s\}$  
is a regular sequence in $A_{i_0}.$ As $\{x_1,\ldots,x_{s-1},x_s\}$ is not a regular sequence in $A_{i_0},$ it cannot be a regular sequence 
in  $A_j,\ \forall j\geq i_0.$ But $x_1,\ldots,x_s$ is a system of parameters in $A_j,\ \forall j\geq i_0.$ We will show that we can assume 
that  $x_2,\ldots,x_s$ remains a regular
sequence in $A_j,\ \forall j\geq i_0.$ If this is so, it follows that $x_2,\ldots,x_s$ is a regular sequence in $A$ and again by
\cite[Theorem 1.7]{K1} we get 
that $A$ is almost Cohen-Macaulay. So, assume that there exists $j_1\geq i_0$ such that $\{x_2,\ldots,x_s\}$ is not a regular sequence
in $A_{j_1}.$ First we observe that then $\{x_2,\ldots,x_s\}$ is not regular in $A_k,\ \forall k\geq j_1.$ Since $A_{j_1}$
is almost Cohen-Macaulay, it follows that, for example, $\{x_1,\widehat{x_2},x_3,\ldots,x_s\}=\{x_1,x_3,\ldots,x_s\}$ is a regular 
sequence in $A_{j_1}.$
Assume that there exists $j_2\geq j_1$ such that $\{x_1,\widehat{x_2},x_3,\ldots,x_s\}$
is not a regular sequence in $A_{j_2}.$ Then there exists $j_2\geq j_1$ such that, for example, $\{x_1,x_2,x_4,\ldots,x_s\}$
is a regular sequence. Continuing we obtain an index $l_0\in\Lambda,$ such that $\{x_1,\ldots,x_{i-1},\widehat{x_i},x_{i+1},\ldots,x_s\}$
is not a regular sequence in $A_l,\ \forall\ l\geq l_0$ and $\forall\ i=1,\ldots,s.$ Contradiction.

\begin{exam}\label{exlim}
Let $K\subset L$ be an infinite algebraic field extension. Then $L=\mathop\bigcup\limits_{i\in\Lambda}K_i,$ where $K_i$ are the finite 
field subextensions of $L.$ For any $i\in\Lambda,$ let $A_i:=K_i[[X^4,X^3Y,XY^3,Y^4]].$ It is easy to see that for any
$i\in\Lambda$ we have  $\de(A_i)=1$ and $\dim(A_i)=2$, that is 
$A_i$ is an almost Cohen-Macaulay ring which is not Cohen-Macaulay. Moreover 
$A:=\mathop\bigcup\limits_{i\in\Lambda}A_i\subsetneq L[[X^4,X^3Y,XY^3,Y^4]]$ and  the family $(A_i)_{i\in\Lambda}$ satisfies 
the conditions in \ref{limite}. 
Hence $A$ is an almost Cohen-Macaulay ring.
\end{exam}
  \par\noindent \textbf{Acknowledgements:}
  This paper was partiallly prepared while the first author was visiting Institut Teknologi Bandung. He thanks this institution for support and hospitality.

\vspace{2.0cm}
\noindent Simion Stoilow Institute of\\
Mathematics of the Romanian\\ Academy\\
P. O. Box 1-764\\
Bucharest 014700\\
Romania\\
email: cristodor.ionescu@imar.ro\\
\vspace{1.0cm}
and\\
Department of Mathematics\\
Payame Noor University\\
P.O. Box 19395-3697\\
Tehran\\
Iran\\
email: samanetabeja\nobreak maat\_golestan@yahoo.com


  \end{document}